\documentstyle[12pt]{article}
\title{Heisenberg Double Description of $\kappa$-Poincar\'{e}
Algebra and $\kappa$-deformed Phase Space \thanks{Supported
by KBN grant 2P30208706}}

\author{J. Lukierski\thanks{Institute for Theoretical Physics,
University 
of Wroc{\l}aw, pl. Maxa Borna 9, 50-204 Wroc{\l}aw, Poland.}
\and A. Nowicki\thanks{ 
Institute of Physics, Pedagogical University,
 pl. S{\l}owia\'{n}ski 6, 65-029 Zielona G\'{o}ra, Poland.}
}
\def\k{\kappa}
\def\poin{Poincar\'e }
\def\kdef{$\k$-deformed }
\date{December 1996}

\begin{document}
\maketitle
\begin{abstract}
The $\kappa$-deformed dual pair of Poincar\'{e} algebra and
Poincar\'{e} group
is formulated in the framework of Heisenberg doubles. The covariant
$\kappa$-deformed  phase space is described in detail as a subalgebra.
The realizations of proposed algebraic scheme are considered.\\
\end{abstract}

\section{The Heisenberg double of $\kappa$-Poincar\'{e} algebra}

It is known that from a pair of dual Hopf algebras describing the
quantum symmetry group ${\cal{A}}$ and its dual quantum Lie algebra
${\cal{A}}^*$ one can construct three different
double algebras: Drinfeld double and Drinfeld codouble, both with 
Hopf 
algebra structure, and Heisenberg double, which is not a Hopf 
algebra. 
Heisenberg double describes the semidirect product of the algebra of
vector
fields acting on the functions of quantum group and represents
algebraic 
generalization of the notion of cotangent double on the group, i.e.
describes 
the generalization of the phase space for the group manifolds.\\

In the present note  we apply the notion of Heisenberg double of a
bialgebra
(see \cite{A}) in the case of the
$\kappa$-deformation of Poincar\'{e} algebra and its dual group
\cite{B,C}. The
Heisenberg double of $\kappa$-Poincar\'{e} algebra $\cal{H(P_\kappa)}$
contains
as its subalgebra another Heisenberg double which turns out to be the
$\kappa$-deformed phase space.\\
We discuss in $\kappa$-deformed phase space the deformed uncertainty
relations
for this phase space which
depend on two parameters: the Planck's constant $\hbar$ and 
the deformation parameter $\k$ with mass dimension. One of the 
important problems in the phase spaceformalism is the canonical change
of phase space coordinates, leaving invariant the Heisenberg canonical
commutation relations. In the Heisenberg double description the object
which 
we would like to leave invariant under canonical transformations are
the scalar products defining the duality relations.
The problem of changing the basis in dual configuration space
introduced by the
nonlinear changes of fourmomentum generators we formulate in the
language of
the finite difference equations.

The deformation of the Heisenberg commutation relations has been
studied
in general algebraic framework by Kempf \cite{kem}. The deformation 
of 
classical relativistic symplectic structure for $\k$-\poin was
discussed 
earlier in \cite{Zak2,E}. Recently the link between the
noncommutativity of 
the spacetime coordinates in \kdef Minkowski space and the effects of 
quantum gravity was considered also by Amelino-Camelia \cite{7}.

\subsection{Heisenberg double}

Let $\{A_a\}$ be a linear basis of bialgebra ${\cal{A}}$ with the
following
\medskip\\
- {\it multiplication relations}:\\
$$
A_a \circ A_b = f_{ab}^{c} A_c
\eqno(1.1a)
$$
- {\it co-multiplication relations}:\\
$$
\Delta(A_a) = h_{a}^{bc} A_b \otimes A_c
\eqno(1.1b)
$$
We recall the following definition \cite{A}\\

{\bf Def.} {\it The Heisenberg double ${\cal{H(A)}}$ of
bialgebra ${\cal{A}}$
is an associative algebra with\\ \hspace*{4em} a basis $\{A_a, B^b\}$
satisfying}\medskip\\ \hspace*{2em}$A_a\circ A_b = f_{ab}^{c}
A_c\hfill(1.2a)$\medskip\\ \hspace*{2em}$B^a\circ B^b = h^{ab}_{c}
B^c$\hfill(1.2b)\medskip\\ \hspace*{2em}$B^a\circ A_b =
f^{a}_{ed}h^{ce}_{b}
A_c \circ B^d\hfill(1.2c)$\medskip\\
The bialgebra ${\cal{A}}$ and dual ${\cal{A}}^*$ are equivalent to
subalgebras
spanned by the basis $\{A_a\}$ and $\{B^b\}$ respectively.\\

The multiplication given by the formula $(1.2c)$ can be derived from
a left
action of a dual algebra ${\cal{A}}^*$ on $\cal{A}$ defined as
follows\\
$$
b\triangleright a = a_{(1)}<b, a_{(2)}> \qquad\qquad b\in {\cal{A}}^*,
a\in
\cal{A} \eqno(1.2d)
$$
where we use Sweedler's notation for coproduct
$\Delta(a)=a_{(1)}\otimes a_{(2)}$.
In particular using duality relation between multiplication and
comultiplication one gets\\
$$
b\triangleright a \tilde{a}
= (a\tilde{a})_{(1)}<b, (a\tilde{a})_{(2)}> =
(a\tilde{a})_{(1)}<\Delta(b), a_{(2)} \otimes {\tilde{a}}_{(2)}> =
$$
$$
=(a\tilde{a})_{(1)}<b_{(1)}, a_{(2)}> <b_{(2)}, {\tilde{a}}_{(2)}> =
a_{(1)} <b_{(1)}, a_{(2)}> ({\tilde{a}}_{(1)}<b_{(2)},
{\tilde{a}}_{(2)}>) =
$$
$$
=a_{(1)}<b_{(1)}, a_{(2)}> (b_{(2)} \triangleright {\tilde{a}}) =
(b\circ a)\triangleright {\tilde{a}}
$$
and we obtain\\
$$
b\circ a = a_{(1)}<b_{(1)}, a_{(2)}> b_{(2)}
\eqno(1.2e)
$$
The multiplication $(1.2c)$ is equivalent to
the commutation relations $(1.2e)$.
In particular, for a linear basis $\{A_a\}$ and $\{B^b\}$
we have\medskip\\
\hspace*{4em}$\Delta(A_a) = h_{a}^{ce} A_{c}\otimes A_{e} = A_{a(1)}
\otimes
A_{a(2)}$\medskip\\
\hspace*{4em}$\Delta(B^b) = f_{cd}^{b} B^{c} \otimes B^{d}
= B_{(1)}^{b} \otimes B_{(2)}^{b}$\medskip\\
and the formula $(1.2e)$ becomes $(1.2c)$.\\

\subsection{$\kappa$-Poincar\'{e} algebra}

We choose the following realization of the $\kappa$-Poincar\'{e}
algebra in bicrossproduct basis 
\cite{B,D}\\$(\mu,\nu,\lambda,\sigma=0,1,2,3; i,j=1,2,3$ 
and $g^{\mu\nu}=g_{\mu\nu}=(-1,1,1,1))$\medskip\\-{\it algebra 
sector}:\medskip\\ \hspace*{2em}$[P_\mu,P_\nu]  =  
0$\medskip\\\hspace*{2em}$[M_{\mu\nu},M_{\lambda\sigma}]
=i\hbar\left(g_{\mu\sigma}M_{\nu\lambda}+
g_{\nu\lambda} M_{\mu\sigma}-g_{\mu\lambda}M_{\nu\sigma}-g_{\nu\sigma}
M_{\mu\lambda}\right)$\medskip\\ \hspace*{2em}$[M_{ij}, P_\mu]  =
-i\hbar\left(g_{i\mu}P_j  - g_{j\mu}P_i\right)$\medskip\\
\hspace*{2em}$[M_{i0}, P_0]  =  i\hbar P_i$\medskip\\
\hspace*{2em}$[M_{i0},
P_j]  = i\delta_{ij}\left({{\hbar}^2}\kappa\sinh({P_{0}\over
\hbar\kappa})e^{-{P_0 \over \hbar\kappa}}
+ {1\over 2\kappa}{\vec{P}}^2 \right)
- {i\over \kappa}P_{i}P_{j}\hfill(1.3)$\medskip\\
-{\it coalgebra sector}:\medskip\\
\hspace*{2em}$\Delta(M_{ij}) \ = \ M_{ij}\otimes I \ + \ I\otimes
M_{ij}$\medskip\\ \hspace*{2em}$\Delta(M_{k0}) \ = \ M_{k0}\otimes
e^{-{P_{0}\over \hbar\kappa}} \ + \ I\otimes M_{k0} \
+ \ {1\over \hbar\kappa}
M_{kl} \otimes P_{l}$\medskip\\ \hspace*{2em}$\Delta(P_{0}) \ =
\ P_{0}\otimes
I \ + \ I\otimes P_{0}$\medskip\\
\hspace*{2em}$\Delta(P_{k}) \ = \ P_{k}\otimes e^{-{P_0\over
\hbar\kappa}}
\ +
\ I\otimes P_k$\hfill(1.4)\medskip\\
where $\hbar$ denotes Plank's constant and $\kappa$ deformation
parameter. The
formulas for the antipode and counit are omitted because they are
not essential
for this construction.

\subsection{$\kappa$-Poincar\'{e} group}

Using the following duality relations\\
$$
<x^{\mu}, P_{\nu}> \ = \ i\hbar \delta^{\mu}_{\nu} \quad
<{\Lambda^{\mu}}_{\nu}, M_{\alpha\beta}> \ = \
i\hbar(\delta^{\mu}_{\alpha}
g_{\nu\beta} \ - \ \delta^{\mu}_\beta g_{\nu\alpha})
\eqno(1.5)
$$
we obtain the commutation relations defining $\kappa$-Poincar\'{e}
group \cite{C,D} in the form\medskip\\
-{\it algebra sector}:\medskip\\
\hspace*{2em}$[x^{\mu}, x^{\nu}] \ = \
{i\over \kappa}(\delta^{\mu}_{0} x^{\nu}
- \delta^{\nu}_{0} x^{\mu})$\medskip\\
\hspace*{2em}$[\Lambda^{\mu}_{\nu}, x^{\lambda}] \ = \ -{i\over
\kappa}\left((\Lambda^{\mu}_{0}
- \delta^{\mu}_{0})\Lambda^{\lambda}_{\nu} +
(\Lambda^{0}_{\nu} - \delta^{0}_{\nu})g^{\mu\lambda}\right)$\medskip\\
\hspace*{2em}$[\Lambda^{\mu}_{\nu}, \Lambda^{\alpha}_{\beta}] \ = \
0\hfill(1.6)$\medskip\\ -{\it coalgebra sector}:\medskip\\
\hspace*{2em}$\Delta(x^{\mu}) \ = \ {\Lambda^{\mu}}_{\alpha}\otimes
x^{\alpha} \ + \ x^\mu\otimes I$\medskip\\
\hspace*{2em}$\Delta({\Lambda^{\mu}}_{\nu}) \ = \
{\Lambda^{\mu}}_{\alpha}\otimes
{\Lambda^{\alpha}}_{\nu}\hfill(1.7)$\medskip\\
The commutation relations $(1.3)$ and $(1.6)$ one can supplement
by the
following relations obtained from (1.2c) or $(1.2e)$\medskip\\
-{\it cross relations}:\medskip\\
\hspace*{2em}$[P_k, x_l] \ = \ -i\hbar\delta_{kl}
\qquad\qquad [P_0, x_0] \ 
= \i\hbar$\medskip\\
\hspace*{2em}$[P_k, x_0] \ = \ -{i\over \kappa} P_k \qquad\qquad
[P_0, x_l] 
\ =\ 0$\medskip\\
\hspace*{2em}$[P_\mu, \Lambda^{\alpha}_{\beta}] \ =
\ 0$\medskip\\
\hspace*{2em}$[M_{\alpha\beta}, \Lambda^{\mu}_{\nu}] \ = \
i\hbar\left(\delta^{\mu}_{\beta} \Lambda_{\alpha\nu}
- \delta^{\mu}_{\alpha}
\Lambda_{\beta\nu}\right)$\medskip\\ \hspace*{2em}$[M_{\alpha\beta},
x^{\mu}] \
= \ i\hbar\left(\delta^{\mu}_{\beta}x_{\alpha} - \delta^{\mu}_{\alpha}
x_{\beta}\right) + {i\over \kappa}\left(\delta^{0}_{\beta}
{M_{\alpha}}^{\mu}
- \delta^{0}_{\alpha}{M_{\beta}}^{\mu}\right)\hfill(1.8)$\medskip\\
where\\
$$
{M_{\alpha}}^{\mu} \ = \ g^{\mu\rho}{M_{\alpha\rho}}\qquad
{M^{\mu}}_{\alpha} \ = \ g^{\mu\rho}{M_{\rho\alpha}}
\eqno(1.9)
$$

The relations $(1.3), (1.6)$ and $(1.8)$ give us the Heisenberg double
$\cal{H(P_\kappa)}$ of $\kappa$-Poincar\'e algebra in terms of
commutation relations.

\section{$\kappa$-Deformed Phase Space}

Let us consider  subalgebra of $\cal{H(P_\kappa)}$ given by the
following\medskip\\
- {\it commutation relations}:\\
$$
\begin{array}{ll}
[x_k , x_l] \ = \ 0\quad &  [P_\mu , P_\nu]  =  0\medskip\\

[x_0 , x_k] \ = \ {i\over \kappa}x_k\quad & \medskip\\

[x_k , P_l] \ = \ i\hbar\delta_{kl}\quad &  [x_k,P_0] \ = \ 
0\medskip\\

[x_0 , P_l] \ = \ {i\over \kappa}P_l\quad &  [x_0 , P_0] \ = \ -i\hbar
\end{array}
\eqno(2.1)
$$
which describe the $\kappa$-deformed phase space.
For $\kappa\rightarrow\infty$
we get the standard nondeformed phase space satisfying the Heisenberg
commutation relations. This phase space one can obtain immediately
using the relations $(1.2)$ to the commuting fourmomentum algebra
with noncocommutative coproduct $(1.4)$.

\subsection{Realization of $\kappa$-deformed phase space algebra
in terms of standard phase space variables}

It is easy to see that the commutation relations $(2.1)$ can
be realized in a Hilbert space in terms of standard undeformed
momentum and position operators satisfying the Heisenberg
commutation relations \cite{E}\medskip\\
$$
[\hat{x}_\mu, \hat{p}_\nu] \ = \ i\hbar g_{\mu\nu}
\eqno(2.2)
$$
then we can define:\medskip\\
\hspace*{4em}$x_k  =  \hat{x}_k \quad  x_0  =  \hat{x}_0 + {1\over
2\kappa\hbar} (\hat{x}_{i}\hat{p}_{i} +
  \hat{p}_{i}\hat{x}_{i})  \quad P_\mu =
\hat{p}_\mu\hfill(2.3)$\medskip\\ Equivalently, we can describe this
transformation introducing the operator\\ $$
U \ = \ e^{{i\over
2\kappa{\hbar^2}}(\hat{x}_{k}\hat{p}_{k}
+\hat{p}_{k}\hat{x}_{k})\hat{p}_0}
\eqno(2.5)
$$
then we get\\
$$
x_k \ = \ \hat{x}_k \qquad x_0 \
= \ U \hat{x}_{0} U^{-1} \qquad P_\mu \
= \ \hat{p}_\mu
\eqno(2.6)
$$
It appears that for an arbitrary invertible operator $U$ , the algebra
generated by $x_\mu, P_\nu$ satisfies the Jacobi identity. Therefore,
such an operator $U$ describes  some general class of deformations of
the standard phase space given by $(2.2)$\\

\subsection{$\kappa$-deformed Heisenberg uncertainty relations}

Introducing  the dispersion of the observable $a$ in quantum
mechanical sense by:\\
$$
\Delta(a) \ = \ \sqrt{<a^2> - <a>^2} \eqno(2.7)
$$
we obtain the nonvanishing $\kappa$-deformed uncertainty relations in 
$\kappa$-deformed phase space as form\medskip\\\hspace*{4em}
$\Delta(x_0)
\Delta(x_k)\geq {1\over2\kappa}|<x_k>|$\medskip\\
\hspace*{4em}$\Delta(P_k) \Delta(x_l)\geq {1\over
2}\hbar\delta_{kl}$\medskip\\
\hspace*{4em}$\Delta(P_0) \Delta(x_0)\geq {1\over
2}\hbar$\medskip\\
\hspace*{4em}$\Delta(P_k) \Delta(x_0)\geq {1\over
2\kappa}|<P_k>|\hfill(2.8)$\medskip\\
which become the standard ones in the limit
$\kappa\rightarrow\infty$.\\

\subsection{Remark on changing the basis}

It is obvious, that for our choice of the basis in phase space,
the following general duality relations hold $(r,s=1,2,3)$:\\
$$
<x^{k}_{0}x^{l}_{s} , P^{n}_{r} P^{m}_{0}> = (\hbar)^{k+l}
\delta_{rs}\delta^{ln}\delta^{km} k!l!(-i)^{k}(i)^{l}
\eqno(2.9a)
$$
$$
<x^{l}_{s} x^{k}_{0}, P^{n}_{r} P^{m}_{0}> =
\left\{\begin{array}{ll}
(\hbar)^{k+l}\delta_{rs}\delta^{ln}{k!l!\over (k-m)!}\left(-{n\over
\kappa\hbar} \right)^{k-m}(-i)^{k}(i)^{l} & \mbox{$ k\geq m$}\\
0 & \mbox{$ k < m$}
\end{array}
\right\}
\eqno(2.9b)
$$
and for any polynomial functions $\psi(x_{0},\vec{x} )$ and
$f(P_{0},\vec{P})$ the duality pairing has the form\medskip\\
$$
<:\psi(x_{0},\vec{x}): , f(P_{0},\vec{P})> =
f(-i\hbar\partial_{0},i\hbar\vec{\bigtriangledown} ) 
\psi(x_{0},\vec{x})|_{x=0}\eqno(2.10)
$$
where $:\psi:$ denotes odered function $\psi$ with all powers
of $x_{0}$ to the left.\\
Let us consider the following change of the basis of
$\kappa$-Poincar\'{e} algebra:\\
$$
\tilde{P}_i = P_{i} e^{P_{0}\over \kappa\hbar} = P_{i} k(P_{0})
\eqno(2.11a)
$$
$$
\tilde{P}_0 = \kappa\hbar\sinh\left(P_{0}\over \kappa\hbar\right)
+ {1\over
2\kappa\hbar}{\vec{P}}^{2}k(P_0) = h(P_0, \vec{P}^2)
\eqno(2.11b)
$$
which transforms to the standard, nondeformed basis of the 
Poincar\'{e}
algebra.\medskip\\
From the duality relations for the phase space, the corresponding
space-time
$\{\tilde{x}_l , \tilde{x}_0\}$ related to the transformed
momentum space
$\{\tilde{P}_l , \tilde{P}_0\}$ i.e.
$$
<\tilde{x}_l , \tilde{P}_k> = i\hbar\delta_{kl}\qquad <\tilde{x}_0 ,
\tilde{P}_0> = -i\hbar
\eqno(2.12)
$$
is given by\\
$$
\tilde{x}_l = :F_{l}(x_0,\vec{x}): \qquad \tilde{x}_0
= :F_{0}(x_0,\vec{x}):
\eqno(2.13)
$$
where the functions $F_l$ and $F_0$ satisfy the following differential
equations:\\
$$
\left({\kappa\over 2}(1-e^{{2i\over \kappa}\partial_0}) - {1\over
2\kappa}\bigtriangleup\right)^{m}\partial^{n}_{r}\left(F_{0}^{k}(x_0
- i{{m+n}\over \kappa}, \vec{x})\cdot\right.
\eqno(2.14)
$$
$$
\left.F_{s}^{l}(x_0 - i{{m+n}\over \kappa},\vec{x})\right)\mid_{x=0} 
= 
(-i)^{k}k!l!\delta_{rs}\delta^{ln}\delta^{km}
\eqno(2.15)
$$
in particular for $l = n = 0$\\
$$
\left({\kappa\over 2}(1-e^{{2i\over \kappa}\partial_0}) - {1\over
2\kappa}\bigtriangleup\right)^{m}F_{0}^{k}(x_0
- i{{m}\over \kappa}, \vec{x})\mid_{x=0} = (-i)^{k}k!\delta^{k m}
\eqno(2.16a)
$$
for $k = m = 0$\\
$$
\partial^{n}_{r}F_{s}^{l}(-i{{n}\over 
\kappa},\vec{x})\mid_{\vec{x}=0}=
 n!\delta_{rs}\delta^{ln}
\eqno(2.16b)
$$
\subsection{Final remarks}

In this note we presented the $\kappa$-deformed phase space for
the standard version of $\kappa$-deformation, which introduces
the "quantized" nature of time coordinate. An analogous discussion
can be presented in the case of
generalized $\kappa$-deformation (see \cite{F,H} and the lecture of 
Maslanka) describing fourdimensional space-time relativistic
symmetries with one arbitrary direction becoming "quantum" after
the $\kappa$-deformation.
It should be recalled that the clasification of the deformations of
D=4 
Poincar\'{e} group has been presented recently by Podles and 
Woronowicz 
\cite{G}. Taking into account the clasification of classical
$r$-matrices 
or D=4 Poincar\'{e} algebra by Zakrzewski \cite{H} it would
be interesting 
to describe all possible deformed quantum phase spaces obtained by 
the 
Heisenberg double construction for known variety of deformed D=4 
fourdimensional relativistic symmetries.\\

\section*{Acknowledgments}

The authors would like to thank the organizers of the 
conference, in particular prof. V.D. Doebner and prof. V.K. Dobrev, 
for 
their hospitality in Goslar.

\end{document}